\documentclass[a4paper]{article}
\usepackage{amsfonts}
\usepackage[fleqn]{amsmath}
\usepackage{amscd}
\usepackage{amssymb}
\usepackage{amsthm}
\usepackage{amsrefs}
\usepackage{latexsym}
\usepackage{bbm}
\usepackage{titlesec}
\usepackage{manfnt}
\numberwithin{equation}{section}
\author{Timothy Foo}
\title{On the parity of the multiplicative order of certain products of integers related to Gauss factorials}
\date{}
\begin{document}
\maketitle
\begin{abstract}
In this note, we prove that under some conditions, certain products of integers related to Gauss factorials are always quadratic residues.
\end{abstract}
\section*{\normalsize Introduction}

Let $n$ be a positive integer and $p$ a prime congruent to $1$ modulo $n$. Then it follows that the numbers $\{1,2,\dots,p-1\}$ may be divided into $(p-1)/n$ equally sized sets. In analogy with Wilson's theorem, it is interesting to ask about the products $\prod_{j=(k-1)(p-1)/n+1}^{k(p-1)/n}j$, where $k=1,\dots,n$. The interesting topic of Gauss factorials was studied in [2-9], and in [5] the partial products to be studied here, which are related to Gauss factorials, were and henceforth here are denoted
$$
\Pi_{k}^{(n)}=\prod_{j=(k-1)(p-1)/n+1}^{k(p-1)/n}j
$$
where $k=1,\dots,n$. In particular, in [2-9], interest has been on the multiplicative order of Gauss factorials, and the accompanying conditions on $p$. In [14], Mordell proved that for $p\equiv 3 \bmod 4$,
$$
\Pi_{1}^{(2)}\equiv (-1)^a \bmod p
$$
where $a\equiv \frac{1}{2}(1+h(-p)) \bmod 2$, $h(-p)$ being the class number of the imaginary quadratic field $\mathbb{Q}[\sqrt{-p}]$. It is interesting to ask, is there a statement for products over various $k$ of $\Pi_{k}^{(q)}$, $q$ an odd prime, that relates to the class number? When $p$ is a prime congruent to $3$ modulo $4$, the multiplicative order of an element in $(\mathbb{Z}/p\mathbb{Z})^{*}$ is either odd or congruent to $2$ modulo $4$. The aim of this note is to show that for primes $p$ congruent to $3$ modulo $4$, and also being congruent to $1$ modulo $q$ so that the appropriate products $\Pi_{k}^{(q)}$, where $k=1,\dots,q$, exist, certain products of these products have odd multiplicative order modulo $p$. For $p$ congruent to $3$ modulo $4$, this is equivalent to these elements being quadratic residues modulo $p$. Thus, we have the following theorem.

\newtheorem*{theorem}{Theorem 1}
\begin{theorem}
Let $q$ be an odd prime. Let $p$ be a prime such that $p\equiv 3 \bmod 4$ and $p\equiv 1 \bmod q$. Then the element
$$
\prod_{\substack{k=1\\ (q+1)/2 -k \equiv 1 \bmod 2}}^{(q-1)/2}\Pi_{k}^{(q)}
$$
is a quadratic residue modulo $p$.
\end{theorem}
\noindent
In the remainder of this note, we shall prove Theorem 1 and discuss its implications. Then we look at further generalizations. Some notation is helpful which we introduce. For $a\in(\mathbb{Z}/p\mathbb{Z})^{*}$, let the Legendre symbol be denoted by
$$
\left(\frac{a}{p}\right) = \begin{cases}
1, \mbox{ if a is a quadratic residue modulo p} \\
-1\mbox{ otherwise.}
\end{cases}
$$
Furthermore, as stated earlier, let $h(-p)$ denote the class number of the imaginary quadratic field $\mathbb{Q}[\sqrt{-p}]$, $p\equiv 3 \bmod 4$. Finally, for $p\equiv 1 \bmod q$, let $a_k$ and $b_k$ denote respectively the number of quadratic residues and nonresidues in the set $\{(k-1)(p-1)/q+1,\dots,k(p-1)/q\}$.

\noindent
At this point, let us mention Theorem 3.1 of [13], an example of which is seen in [12]. Specifically, it follows from this theorem that if $q\equiv 3 \bmod 4$, $H$ the subgroup of squares in $(\mathbb{Z}/q\mathbb{Z})^{*}$, $\beta = (\sum_{\substack{i=1\\-i\in H}}^{q-1}i)/q$, $p\equiv 1 \bmod q$, $4p^{h(-q)}=a^2+qb^2$, where $a\equiv 2 \bmod q$, then 
\begin{equation}\label{a}
\left(\frac{a}{p}\right) = \left(\frac{\left(-1\right)^{\beta} \prod_{-i\in H}\left(i\left(\frac{p-1}{q}\right)\right)!}{p}\right).
\end{equation}
We have that
\newtheorem*{theorem2}{Theorem 2}
\begin{theorem2}
Let $p,q$ be primes, $q\equiv 3 \bmod 4$, $p\equiv 1 \bmod q$, $p\equiv 3\bmod 4$. Let $a$ be defined as above. Then the following are equivalent.\newline
(1)$$\mbox{Theorem 1 applied to the case when q is 3 mod 4} \Leftrightarrow \mbox{Equation }(\ref{a})$$
\newline
(2)$$ \left(\frac{a}{p}\right) = (-1)^{(q+1)/4}.$$
\end{theorem2}
{\it Proof}. 
Since $\beta=\frac{h(-q)+1}{2}+\frac{q-3}{4}$ and $q\equiv 3\bmod 4$,
\begin{eqnarray*}
1&=&\left(\frac{\prod_{\substack{k=1\\ (q+1)/2 -k \equiv 1 \bmod 2}}^{(q-1)/2}\Pi_{k}^{(q)}}{p}\right)\mbox{ by Theorem 1}\\
&=&\left(\frac{\prod_{i=1}^{(q-1)/2}\left(i\left(\frac{p-1}{q}\right)\right)!}{p}\right)\\
&=&(-1)^{1+(h(-q)+1)/2}\left(\frac{\prod_{-i\in H}\left(i\left(\frac{p-1}{q}\right)\right)!}{p}\right)\\
&=&(-1)^{(q+1)/4}\left(\frac{a}{p}\right) \mbox{ by Equation }(\ref{a}).
\end{eqnarray*}

\section*{\normalsize The expression for the class number of an imaginary quadratic field in terms of Legendre symbols}

The proof of Theorem 1 here begins with the beautiful expression of the class number of the imaginary quadratic field $\mathbb{Q}[\sqrt{-p}]$, $p\equiv 3 \bmod 4$, due to Dirichlet. The formula is given by
$$
h(-p) = \frac{1}{2-\left(\frac{2}{p}\right)}\sum_{a=1}^{(p-1)/2}\left(\frac{a}{p}\right)
$$
where $p\equiv 3 \bmod 4$, $p>3$ and references include Corollary 5.3.13 of [1], and [10] (equations (7), (8) of Chapter 1 and equation (19) of chapter 6). This nice formula is due to the Legendre symbol being an odd character when $p\equiv 3\bmod 4$. A different form of this formula, generalized to the other primes, instead of the prime $2$, may be obtained by similar treatment of the L-function $L(1,(\frac{}{p}))$ following the outline summarized in Chapter 1 of [10]. This becomes our first lemma.

\newtheorem*{lemma}{Lemma 1}
\begin{lemma}
Let $q$ be an odd prime, and $p\equiv 3\bmod 4$, $p>3$. Then for $p\not=q$,
$$
h(-p) = \frac{1}{q-\left(\frac{q}{p}\right)}\sum_{a=1}^{(p-1)/2}\left(\frac{a}{p}\right)\left(q-1-2\left\lfloor\frac{aq}{p}\right\rfloor\right).
$$
\end{lemma}

\section*{\normalsize Proof of Theorem 1}

Before proceeding further, we recollect the fact that the class number $h(-p)$ is odd, which is our second lemma.

\newtheorem*{lemma2}{Lemma 2}
\begin{lemma2}
The class number $h(-p)$ of the imaginary quadratic field $\mathbb{Q}[\sqrt{-p}]$, $p\equiv 3 \bmod 4$, is odd.
\end{lemma2}
\noindent
{\it Proof.} This fact is stated in [11].
\newline
Proceeding to the proof of the theorem, we find that when $p\equiv 1 \bmod q$, Lemma 1 may be written as 
\begin{equation}\label{1}
\left(q-\left(\frac{q}{p}\right)\right)h(-p)/2 = \sum_{k=1}^{(q-1)/2}(a_k-b_k)\left(\frac{q+1}{2}-k\right).
\end{equation}
Noting that $a_k + b_k=\frac{p-1}{q}$, $k=1,\dots,(q-1)/2$, we may combine this with equation \ref{1} to obtain

\begin{equation}\label{2}
\left(\frac{q^2-1}{8}\right)\left(\frac{p-1}{2q}\right) - \frac{\left(q-\left(\frac{q}{p}\right)\right)h(-p)}{4} = \sum_{k=1}^{(q-1)/2}b_k\left(\frac{q+1}{2}-k\right).
\end{equation}
\noindent
The left hand side of equation \ref{2} is always even.\\
\newline
\noindent
{\bf Case $q\equiv 1\bmod 8$}. Then since $p\equiv 1\bmod q$, by quadratic reciprocity, $\frac{\left(q-\left(\frac{q}{p}\right)\right)h(-p)}{4}$ is even, and considering $q\equiv 1,9$ mod $16$ shows that $\left(\frac{q^2-1}{8}\right)\left(\frac{p-1}{2q}\right)$ is even.
\newline
\noindent
{\bf Case $q\equiv 5\bmod 8$}. Then since $p\equiv 1\bmod q$, by quadratic reciprocity, and by Lemma 2, $\frac{\left(q-\left(\frac{q}{p}\right)\right)h(-p)}{4}$ is odd. Considering $q\equiv 5,13$ mod $16$ shows that $\left(\frac{q^2-1}{8}\right)\left(\frac{p-1}{2q}\right)$ is also odd. 
\newline
\noindent
{\bf Case $q\equiv 3\bmod 8$}. Then since $p\equiv 1\bmod q$, by quadratic reciprocity, and by Lemma 2, $\frac{\left(q-\left(\frac{q}{p}\right)\right)h(-p)}{4}$ is odd. Considering $q\equiv 3,11$ mod $16$ shows that $\left(\frac{q^2-1}{8}\right)\left(\frac{p-1}{2q}\right)$ is also odd.  
\newline
\noindent
{\bf Case $q\equiv 7\bmod 8$}. Then since $p\equiv 1\bmod q$, by quadratic reciprocity, $\frac{\left(q-\left(\frac{q}{p}\right)\right)h(-p)}{4}$ is even. Considering $q\equiv 7,15$ mod $16$ shows that $\left(\frac{q^2-1}{8}\right)\left(\frac{p-1}{2q}\right)$ is even.
\newline
\noindent
Reducing equation \ref{2} modulo $2$, we obtain 
$$
\sum_{\substack{k=1\\(q+1)/2 -k \equiv 1 \bmod 2}}^{(q-1)/2}b_k\equiv 0\bmod 2.
$$
This proves Theorem 1.

\section*{\normalsize Relevance for multiplicative orders of partial products and applications}

In general, Wilson's theorem tells us that the product of all the products
$$
\prod_{k=1}^{q}\Pi_{k}^{(q)}
$$
is a quadratic non-residue for $p\equiv 3\bmod 4$. Let $q$ be an odd prime. Let $p\equiv 3 \bmod 4$ and $p\equiv 1 \bmod q$. Since $(p-1)/q$ is even, one has $\Pi_{k}^{(q)}\equiv\Pi_{q-k}^{(q)}\bmod p$, so $\Pi_{(q+1)/2}^{(q)}$ is a non-residue. It is thus interesting to ask which of the products $\Pi_{k}^{(q)}$, $k=1,\dots,(q-1)/2$ are quadratic residues, in which case they will have odd multiplicative order modulo $p$. We have shown that for $p$ satisfying $p\equiv 3\bmod 4$ and $p\equiv 1 \bmod q$, the product of the products $\prod_{\substack{k=1\\ (q+1)/2 -k \equiv 1 \bmod 2}}^{(q-1)/2}\Pi_{k}^{(q)}$ always has odd multiplicative order modulo $p$. Moreover, we may say the following.

\newtheorem*{corollary}{Corollary}
\begin{corollary}
Let $q$ be an odd prime. Let $p\equiv 3 \bmod 4$ and $p\equiv 1 \bmod q$. Consider the set of products $\{\Pi_{k}^{(q)}: k=1,\dots,(q-1)/2, \mbox{ and }(q+1)/2-k\equiv 1 \bmod 2\}$. The number of elements in this set that are quadratic non-residues is even. Thus, the number of elements in this set that have even multiplicative order modulo $p$ is even.
\end{corollary}
\noindent
As examples of the corollary, we consider some cases for small $q$.
\newline
\noindent
{\bf Case $q=3$}. Then $\Pi_{1}^{(3)}$ is a quadratic residue and has odd order modulo $p$.
\newline
\noindent
{\bf Case $q=5$}. Then $\Pi_{2}^{(5)}$ is a quadratic residue and has odd order modulo $p$.
\newline
\noindent
{\bf Case $q=7$}. Then $\Pi_{1}^{(7)}\Pi_{3}^{(7)}$ is a quadratic residue and has odd order modulo $p$. Either both $\Pi_{1}^{(7)}$, $\Pi_{3}^{(7)}$ are quadratic residues, or both are quadratic non-residues.\newline\noindent
{\bf Case $q=11$}. Then $\Pi_{1}^{(11)}\Pi_{3}^{(11)}\Pi_{5}^{(11)}$ is a quadratic residue modulo $p$.

\section*{\normalsize Further generalizations}

Theorem 1 may be generalized to primes $p$ such that $p\equiv 3 \bmod 4$ and $p\not\equiv 1 \bmod q$. In this case, we need a new generalization of the partial products $\Pi_{k}^{(q)}$. With this in mind, set 
$$
\Pi^{\prime (q)}_{k} = \prod_{j=\lfloor\frac{(k-1)p}{q}\rfloor+1}^{\lfloor\frac{kp}{q}\rfloor}j,
$$
$1\leq k \leq q-1$,
$$
\Pi^{\prime (q)}_{q} = \prod_{j=\lfloor\frac{(q-1)p}{q}\rfloor+1}^{p-1}j.
$$
We have that $\Pi^{\prime (q)}_{k}=\Pi_{k}^{(q)}$ when $p\equiv 1 \bmod q$. Let $a_k^{\prime}, b_k^{\prime}$, $1\leq k\leq q-1$ denote respectively the number of quadratic residues and nonresidues in $\{\lfloor\frac{(k-1)p}{q}\rfloor+1, \dots,\lfloor\frac{kp}{q}\rfloor\}$, and $a_q^{\prime}, b_q^{\prime}$ denote respectively the number of quadratic residues and nonresidues in $\{\lfloor\frac{(q-1)p}{q}\rfloor+1,\dots,p-1\}$. These agree with the definition of $a_k,b_k$ when $p\equiv 1 \bmod q$. Let us consider the case $p\equiv 2 \bmod q$. Then we have

\newtheorem*{theorem3}{Theorem 3}
\begin{theorem3}
Let $q$ be an odd prime. Let $p$ be a prime such that $p\equiv 3 \bmod 4$ and $p\equiv 2 \bmod q$. Then the element
$$
b=\prod_{\substack{k=1\\ (q+1)/2 -k \equiv 1 \bmod 2}}^{(q-1)/2}\Pi_{k}^{\prime (q)}
$$
satisfies
$$
\left(\frac{b}{p}\right) = 
\begin{cases}
1, \mbox{ if }q\equiv \pm 1 \bmod 16,\\
-1, \mbox{ if }q\equiv \pm 7 \bmod 16\\
(-1)^{\frac{h(-p)+1}{2}}, \mbox{ if }q\equiv \pm 3 \bmod 16\\
(-1)^{1+\frac{h(-p)+1}{2}}, \mbox{ if }q\equiv \pm 5 \bmod 16.
\end{cases}
$$
\end{theorem3}
\noindent
To verify this, we begin by checking that $a_k^{\prime}+b_k^{\prime}=\frac{p-2}{q}$, $1\leq k \leq q$ and $k\not= (q+1)/2$, and $a_k^{\prime}+b_k^{\prime}=\frac{p-2}{q}+1$ when $k=(q+1)/2$.  Combining with Lemma 1 yields
$$
\left(\frac{q^2-1}{8}\right)\left(\frac{p-2}{2q}\right) - \frac{\left(q-\left(\frac{q}{p}\right)\right)h(-p)}{4} = \sum_{k=1}^{(q-1)/2}b_k^{\prime}\left(\frac{q+1}{2}-k\right).
$$
\noindent
{\bf Case $q\equiv 1\bmod 16$}. Then $\frac{\left(q-\left(\frac{q}{p}\right)\right)h(-p)}{4}$ is even, and checking $q\bmod 32$ shows that $\left(\frac{q^2-1}{8}\right)\left(\frac{p-2}{2q}\right)$ is even.
\newline
{\bf Case $q\equiv 9\bmod 16$}. Then $\frac{\left(q-\left(\frac{q}{p}\right)\right)h(-p)}{4}$ is even, and checking $q\bmod 32$ shows that $\left(\frac{q^2-1}{8}\right)\left(\frac{p-2}{2q}\right)$ is odd.\newline
{\bf Case $q\equiv 7\bmod 16$}. Then $\frac{\left(q-\left(\frac{q}{p}\right)\right)h(-p)}{4}$ is even, and checking $q\bmod 32$ shows that $\left(\frac{q^2-1}{8}\right)\left(\frac{p-2}{2q}\right)$ is odd.\newline
{\bf Case $q\equiv 15\bmod 16$}. Then $\frac{\left(q-\left(\frac{q}{p}\right)\right)h(-p)}{4}$ is even, and checking $q\bmod 32$ shows that $\left(\frac{q^2-1}{8}\right)\left(\frac{p-2}{2q}\right)$ is even.\newline
{\bf Case $q\equiv 5\bmod 16$}. Then  $\frac{\left(q-\left(\frac{q}{p}\right)\right)h(-p)}{4}$ is a half-integer with numerator congruent to $3h(-q)\bmod 4$. Checking $q\bmod 32$ shows that $\left(\frac{q^2-1}{8}\right)\left(\frac{p-2}{2q}\right)$ is a half-integer with numerator congruent to $3\bmod 4$. \newline
{\bf Case $q\equiv 13\bmod 16$}. Then  $\frac{\left(q-\left(\frac{q}{p}\right)\right)h(-p)}{4}$ is a half-integer with numerator congruent to $3h(-q)\bmod 4$. Checking $q\bmod 32$ shows that $\left(\frac{q^2-1}{8}\right)\left(\frac{p-2}{2q}\right)$ is a half-integer with numerator congruent to $1\bmod 4$. \newline
{\bf Case $q\equiv 3\bmod 16$}. Then  $\frac{\left(q-\left(\frac{q}{p}\right)\right)h(-p)}{4}$ is a half-integer with numerator congruent to $h(-q)\bmod 4$. Checking $q\bmod 32$ shows that $\left(\frac{q^2-1}{8}\right)\left(\frac{p-2}{2q}\right)$ is a half-integer with numerator congruent to $3\bmod 4$. \newline
{\bf Case $q\equiv 11\bmod 16$}. Then  $\frac{\left(q-\left(\frac{q}{p}\right)\right)h(-p)}{4}$ is a half-integer with numerator congruent to $h(-q)\bmod 4$. Checking $q\bmod 32$ shows that $\left(\frac{q^2-1}{8}\right)\left(\frac{p-2}{2q}\right)$ is a half-integer with numerator congruent to $1\bmod 4$. 
\newline
\noindent
Next is the case $p\equiv 3 \bmod q$.

\newtheorem*{theorem4}{Theorem 4}
\begin{theorem4}
Let $q>3$ be an odd prime. Let $p$ be a prime such that $p\equiv 3 \bmod 4$ and $p\equiv 3 \bmod q$. Then the element
$$
b=\prod_{\substack{k=1\\ (q+1)/2 -k \equiv 1 \bmod 2}}^{(q-1)/2}\Pi_{k}^{\prime (q)}
$$
satisfies
$$
\left(\frac{b}{p}\right) = 
\begin{cases}
1, \mbox{ if }q\equiv \pm 1 \bmod 12,\\
(-1)^{\frac{h(-p)+1}{2}}, \mbox{ if }q\equiv \pm 5 \bmod 12\\
\end{cases}
$$
\end{theorem4}
\noindent
To verify this, we begin by checking that $a_k^{\prime}+b_k^{\prime}=\frac{p-3}{q}$, $1\leq k \leq (q-1)/2$ and $k\not= \frac{q+2+\left(\left(\frac{q}{3}\right)-1\right)/2}{3}$, and $a_k^{\prime}+b_k^{\prime}=\frac{p-3}{q}+1$ when $k=\frac{q+2+\left(\left(\frac{q}{3}\right)-1\right)/2}{3}$.  Combining with Lemma 1 yields

$$
\left(\frac{q^2-1}{8}\right)\left(\frac{p-3}{2q}\right) + \frac{q-\left(\frac{q}{3}\right)}{12}- \frac{\left(q-\left(\frac{q}{p}\right)\right)h(-p)}{4} = \sum_{k=1}^{(q-1)/2}b_k^{\prime}\left(\frac{q+1}{2}-k\right),
$$
the left hand side of the above being
$$
\frac{\left(q-\left(\frac{q}{3}\right)\right)(1-3h(-p))}{12}
$$
mod $2$.

\section*{\normalsize References}
[1] H. Cohen, A Course in Computational Algebraic Number Theory, Springer Graduate Texts in Mathematics, No. 138, Berlin-Heidelberg 1993.\newline
[2] J. B. Cosgrave and K. Dilcher, Extensions of the Gauss-Wilson Theorem, Integers: Electronic Journal of Combinatorial Number Theory 8 (2008), A39.\newline
[3] J. B. Cosgrave and K. Dilcher, Mod $p^3$ Analogues of Theorems of Gauss and Jacobi on Binomial Coefficients, Acta Arithmetica, Vol. 142, No. 2, 103--108, 2010.\newline
[4] J. B. Cosgrave and K. Dilcher, The Multiplicative Order of Certain Gauss Factorials, International Journal of Number Theory, Vol. 7 (2011), 145--171.\newline
[5] J. B. Cosgrave and K. Dilcher, An Introduction to Gauss Factorials, The American Mathematical Monthly, Vol. 118, No. 9, 812--829.\newline
[6] J. B. Cosgrave and K. Dilcher, Sums of Reciprocals Modulo Composite Integers, Journal of Number Theory, 133, (2013), 3565--3577.\newline
[7] J. B. Cosgrave and K. Dilcher, The Gauss-Wilson Theorem for Quarter-Intervals, Acta Mathematica Hungarica, 142, (2014), 199--230.\newline
[8] J. B. Cosgrave and K. Dilcher, The Multiplicative Order of Certain Gauss Factorials II, Functiones et Approximatio Commentarii Mathematici.\newline
[9] J. B. Cosgrave and K. Dilcher,  A Role for Generalized Fermat Numbers, Mathematics of Computation, September 2015.\newline
[10] H. Davenport, Multiplicative Number Theory, Second Edition, Graduate Texts in Mathematics, Vol. 74, Springer-Verlag, New York, 1980.\newline
[11] A. Laradji, M. Mignotte, and N. Tzanakis, Elementary Trigonometric Sums Related to Quadratic Residues, Elemente Der Mathematik, Vol. 67, Issue 2, 2012, 51--60.\newline
[12] S. G. Hahn and D. H. Lee, Some Congruences for Binomial Coefficients, Class Field Theory - Its Centenary and Prospect, Advanced Studies in Pure Mathemmatics, 30, Math. Soc. Japan (2001).\newline
[13]  S. G. Hahn and D. H. Lee, Gauss Sums and Binomial Coefficients, Journal of Number Theory, Volume 92, Issue 2, 2002, 257--271.\newline
[14] L. J. Mordell, The Congruence $((p-1)/2)!\equiv \pm 1 (\bmod p)$, American Mathematical Monthly 68 (1961), 145--146.
\vspace*{.7cm}

\noindent\begin{tabular}{p{8cm}p{8cm}}
School of Computing, &\\
National University of Singapore, Singapore 119077 \\
Email: {\tt tchfoo@hotmail.com} \\
\end{tabular}
\end{document}